\documentclass{article}

\usepackage{amsthm}
\usepackage{amsmath}
\usepackage{amssymb}

\newtheorem{theorem}{Theorem}
\newtheorem{lemma}[theorem]{Lemma}

\newtheorem{definition}{Definition}

\newtheoremstyle{example}{\topsep}{\topsep}%
{}
{}
{\bfseries}
{.}
{2mm}
{\thmname{#1}\thmnumber{ #2}\thmnote{ #3}}

\theoremstyle{example}
\newtheorem{example}[theorem]{Example}

\newcommand{\R}{\mathbb{R}}
\newcommand{\C}{\mathbb{C}}
\newcommand{\D}{\mathbb{D}}
\renewcommand{\Re}{\mathrm{Re}}
\newcommand{\dbar}{\overline{\partial}}

\title{Variational Formulas for the Green Function}
\date{\today}

\author{Charles Z. Martin}

\begin{document}

\maketitle

\begin{abstract}
  The Green function has a complex dependence upon its underlying domain and differential operator. We briefly review Hadamard's formula
  for the first variation of the Green function due to a perturbation of the domain. We then take a different avenue and approximate the change
  in the Green function when the Laplacian is perturbed into a number of different operators: Helmholtz, Schr\"odinger and Laplace--Beltrami.
\end{abstract}

\maketitle
\section{Introduction}

The study of boundary value problems presents a significant computational challenge. In theory there are numerous formulas, transforms and methods for solving, say, the
Dirichlet problem on an ellipse for the Laplace--Beltrami operator. In fact, knowledge of a single function---the Green function---would reduce the
problem to that of calculating an integral. However, the computation of a Green function is prohibitively difficult; considering its central role
we might turn to methods of approximating it instead. One approach aims to replace the Green function of the problem
with that of an easier problem which is, in some sense, nearby. There are two ways this heuristic can proceed: variation of the underlying domain and variation of the relevant
operator.

The purpose of domain variation is clear; the Green function of a region is nearly impossible to compute for all but the nicest of domains. Approximating a given domain
by a simpler one is a natural way to circumvent the difficulty. The resulting alteration to the Green function is very complicated, but a formula due to Hadamard gives
the first variation, which can be sufficient information for some applications.

For other applications, we might want to approximate the Green function not under perturbations of the domain, but rather of the operator. This avenue of thought is less studied
and for good reason; in many applications the differential operator is considered to be well-known. If we turn our attention to inverse problems however, we might
pose questions wherein an unknown operator has a Green function with some given property. The existence of the underlying operator can be studied by understanding
what sorts of variations in the Green function can result from perturbing a well--understood operator, such as the Laplacian.

In this article we briefly review a formal argument developing Hadamard's variation. With this in hand
we can make a few remarks about what sorts of basic facts can follow from such a formula.
From here we proceed to consider Helmholtz, Schr\"odinger and Laplace--Beltrami operators which are somehow close to the Laplacian. The resulting variation formulas
are new---or at least not widely known---and have nice forms open to interpretation, both mathematical and physical.

\section{Definitions and Known Results}

There are numerous definitions and conventions surrounding the Green function of a domain. For clarity we take a moment to review a few important properties which will
be used in the following.
Let $D\subset \C$ be a bounded domain. A Green function $g:D\times D\to [-\infty,\infty)$ is a map satisfying three properties:
\begin{enumerate}
  \item For every $w\in D$, $g(\cdot,w)$ is harmonic on $D\setminus\{w\}$ and bounded outside every neighborhood of $w$.
  \item For every $w\in D$ we have $g(w,w) = -\infty$ and $g(z,w) = (2\pi)^{-1}\log|z-w| + O(1)$ as $z\to w$.
  \item For any $w\in D$ we have $\lim\limits_{z\to \zeta}g(z,w) = 0$ for quasi-every $\zeta\in \partial D$; that is, for all $\zeta$ except a set of capacity zero.
\end{enumerate}
The Green function exists and is unique for any domain whose boundary has positive capacity. In particular, we will assume that $D$ has real analytic boundary.
Occasionally we wish to emphasize the dependence on one argument; hence we sometimes write $g_w(z)$ in place of $g(z,w)$.
We have the following facts:
\begin{enumerate}
  \item For any $\zeta \in \partial D$ which possesses a barrier, $g(z,w) = 0$ as $z\to\zeta$; having assumed $\partial D$ real analytic, this follows for all $\zeta$.
  We can thus extend $g$ to allow one of its arguments to be on $\partial D$.
  
  \item The Green function is the solution to the boundary value problem
  \begin{equation*}
    \begin{cases}
      \Delta g_w = \delta_w & \textrm{ in } D\\
      g_w = 0 & \textrm{ on } \partial D
    \end{cases}
  \end{equation*}
  It follows that the solution to a more general problem
  \begin{equation*}
    \begin{cases}
      \Delta u = f & \textrm{ in } D\\
      u = 0 & \textrm{ on } \partial D
    \end{cases}
  \end{equation*}
  is given by integration against the Green function:
  \begin{equation*}
    u(z) = \int_D g_z f\, dA.
  \end{equation*}
  This property can be taken as the definition of the Green function of $D$. We can define the Green function $g^*$ of an operator $L$ analogously:
  \begin{equation*}
    \begin{cases}
      L g^*_w = \delta_w & \textrm{ in } D\\
      g^*_w = 0 & \textrm{ on } \partial D
    \end{cases}
  \end{equation*}
  
  \item The Green function is subharmonic.
  \item For all $z,w\in D$ the Green function is symmetric: $g(z,w)=g(w,z)$.
  \item Along $\partial D$ we have
  \begin{equation*}
    \frac{\partial g}{\partial n}(\zeta,w) \,ds(\zeta)= d\omega(w,\zeta)
  \end{equation*}
  where $\omega$ represents the harmonic measure of $D$ and $\partial g/\partial n$ represents the $\zeta$ derivative in the outward normal direction.
\end{enumerate}
In what follows, we will need the Poisson-Jensen formula for subharmonic functions. Suppose $u$ is subharmonic in a neighborhood of $\overline{D}$ and $u \not\equiv -\infty$. Then
for all $z\in D$
\begin{equation*}
  u(z) = \int_{\partial D} u(\zeta)\, d\omega(z,\zeta) + \int_D g(z,\xi) \Delta u(\xi)\, dA(\xi).
\end{equation*}
Even if $\Delta u$ is a distribution, the formula still holds by interpreting the area integral as a distributional pairing.
Essentially this result follows from the Green identity
\begin{equation*}
  \int_{\partial D} \left(u\frac{\partial v}{\partial n} - v\frac{\partial u}{\partial n}\right)\, ds = \int_{D} \left(u\Delta v - v\Delta u\right)\, dA
\end{equation*}
with $v=g_z$.

\section{Domain Variation of the Green Function}

Before proceeding to variations of operators, we examine Hadamard's formula for the variation of the Green function due to change of the domain.
Suppose that $D$ is slightly enlarged into a new domain $D^*$ in such a way so that along $\partial D$ each point is moved along the outward normal direction a small distance.
The following first--order variation formula gives a relationship between the new Green function $g^*$ and the old one. Because this is well--known, we only give a formal argument;
the original derivation can be found in \cite{Hadamard} and a more rigorous treatment can be found in \cite{SS} or in chapter 15 of \cite{Garabedian}.

\begin{theorem}[Hadamard's formula] \label{Hadamard}
  Let $p\in C(\partial D)$ be a positive function and suppose that for $\epsilon > 0$ each point $\zeta\in\partial D$ is moved along the outward normal direction
  a distance $\epsilon p(\zeta)$. The Green function $g^*$ of the new domain $D^*$ satisfies
  \begin{equation*}
    g^*(z,w) - g(z,w) = -\epsilon\int_{\partial D} \frac{\partial g}{\partial n}(z,\zeta)\frac{\partial g}{\partial n}(w,\zeta) p(\zeta)\,ds(\zeta) + o(\epsilon).
  \end{equation*}
\end{theorem}

\begin{proof}
Fix $w\in D$ and apply the Poisson--Jensen formula to $g^*_w$:
  \begin{equation*}
    g^*(z,w) = \int_{\partial D} g^*(\zeta,w)\, d\omega(z,\zeta) + \int_D g(z,\xi) \Delta g^*(\xi,w)\, dA(\xi).
  \end{equation*}
  Since $\Delta g^*_w = \delta_w$ and $g^*$ is symmetric in its arguments, this becomes
  \begin{equation}
    g^*(z,w) = \int_{\partial D} g^*(w, \zeta)\, d\omega(z,\zeta) + g(z,w). \label{varpf}
  \end{equation}
  Denote the outward unit normal by $n$. For $\zeta \in \partial D$, $\zeta + \epsilon p(\zeta)n \in \partial D^*$. Therefore
  \begin{equation*}
    g^*(w, \zeta  + \epsilon p(\zeta)n) = 0.
  \end{equation*}
  From here we obtain the first-order approximation
  \begin{equation*}
    g^*(w,\zeta) + \epsilon \frac{\partial g^*}{\partial n}(w,\zeta)  p(\zeta) + o(\epsilon) = 0.
  \end{equation*}
  Using this approximation \eqref{varpf} becomes 
  \begin{equation}
    g^*(z,w) - g(z,w) = -\epsilon\int_{\partial D} \frac{\partial g}{\partial n}(z,\zeta)\frac{\partial g^*}{\partial n}(w,\zeta) p(\zeta)\,ds(\zeta) + o(\epsilon), \label{varpf2}
  \end{equation}
  where we have substituted $\partial g/\partial n(z,\zeta) \, ds(\zeta)$ for $d\omega(z,\zeta)$.
  Finally, we remove the dependence on $g^*$ by iterating the formula; that is, from \eqref{varpf2} we have $g^*-g = o(1)$ and
  \begin{equation*}
    g^*(z,w) -g(z,w) = -\epsilon\int_{\partial D} \frac{\partial g}{\partial n}(z,\zeta)\frac{\partial}{\partial n}\left[g(w,\zeta) + o(1)\right](w,\zeta) p(\zeta)\,ds(\zeta) + o(\epsilon).
  \end{equation*}
  This gives the result.
\end{proof}

In traditional notation of the calculus of variations, we can write the aforementioned formula as
\begin{equation*}
	\delta g(z,w) :=\lim_{\epsilon\to 0}\frac{g^*(z,w)-g(z,w)}{\epsilon}
	= -\int_{\partial D} \frac{\partial g}{\partial n}(z,\zeta)\frac{\partial g}{\partial n}(w,\zeta) p(\zeta)\,ds(\zeta).
\end{equation*}
We will at times use this notation for sake of clarity.

\begin{example}
Consider the unit disk $\D\subset \C$. The Green function for a disk of radius $R$ is
\begin{equation*}
  g(z,w) = \frac{1}{2\pi}\ln\left| \frac{R(z-w)}{R^2 - z\bar{w}} \right|.
\end{equation*}
Suppose that $\D$ is enlarged by uniformly increasing the radius by $\epsilon$---that is, the perturbation function $p$ is constant along $\partial D$.
On one hand, a direct computation shows
\begin{equation*}
  \delta g(z,w) = \frac{\partial g}{\partial R}(z,w)\bigg|_{R=1} = -\frac{1}{2\pi} \frac{1-|zw|^2}{ |1-z\bar{w}|^2}.
\end{equation*}
The variational formula gives
\begin{align*}
  \delta g(z,w) &= -\frac{1}{2\pi}\int_{\partial \D} \frac{(1-|w|^2)(1-|z|^2)}{|\zeta-w|^2|\zeta-z|^2}  \, ds(\zeta) \\
  &= -\frac{1}{2\pi}\int_{\partial \D} \frac{(1-|w|^2)(1-|z|^2)}{(\zeta-w)(\zeta-z)(\bar{\zeta} - \bar{z})(\bar{\zeta} - \bar{w})}  \, \frac{d\zeta}{i\zeta}.
\end{align*}
Since $|\zeta| = 1$ we can rewrite $\bar{\zeta} - \bar{z} = (1-\zeta \bar{z})/\zeta$. The above integral can be written as
\begin{equation*}
  \delta g(z,w) = -\frac{1}{2\pi i}\int_{\partial \D} \frac{1}{(\zeta-w)(\zeta-z)} \cdot \frac{(1-|w|^2)(1-|z|^2)\zeta}{(1 - \zeta\bar{z})(1 - \zeta\bar{w})}  \, d\zeta.
\end{equation*}
Note that the second factor in the integral is analytic in a neighborhood of $\overline{\D}$. The Cauchy integral theorem gives
\begin{align*}
  \delta g(z,w) &=-\frac{1}{2\pi i}\int_{\partial D} \frac{1}{z-w}\left(\frac{1}{\zeta-z} - \frac{1}{\zeta-w}\right)
  \cdot \frac{(1-|w|^2)(1-|z|^2)\zeta}{(1 - \zeta\bar{z})(1 - \zeta\bar{w})}  \, d\zeta +o(\epsilon) \\
  &= -\frac{1}{z-w}\left( \frac{z(1-|w|^2)}{1 - z\bar{w}} - \frac{w(1-|z|^2)}{1 - w\bar{z}} \right) \\
  &= -\frac{1}{2\pi} \frac{1-|zw|^2}{ |1-z\bar{w}|^2},
\end{align*}
so the two answers agree.
\end{example}

Before proceeding we make a few remarks.
\begin{enumerate}
  \item The outward normal derivative of the Green function is positive since, for instance, it is the density of the domain's harmonic measure.
  Therefore the variation is always negative, so we conclude that enlarging a domain decreases the Green function at every point.
  
  \item Let $\lambda :\C\to\R$ be a positive smooth function and define the operator $L = \nabla\lambda\nabla$. If we alter the definition of the Green function so that
  $Lg = \delta$, then we can derive another variational formula. From Green's identity
  \begin{equation*}
    \int_{\partial D} \lambda\left(u\frac{\partial v}{\partial n} - v\frac{\partial u}{\partial n}\right)\, ds = \int_{D} \left(uL v - vL u\right)\, dA,
  \end{equation*}
  we set $v = g_z$ to get an analogue of the Poisson-Jensen formula:
  \begin{equation*}
    u(z) = \int_{\partial D} u(\zeta) \lambda(\zeta) \, d\omega(z,\zeta) + \int_D g(z,\xi) \Delta u(\xi)\, dA(\xi).
  \end{equation*}
  From here we obtain the first variation of $g$.
  \begin{equation*}
    g^*(z,w) -  g(z,w) = -\epsilon\int_{\partial D} \frac{\partial g}{\partial n}(z,\zeta)\frac{\partial g}{\partial n}(w,\zeta) \lambda(\zeta)p(\zeta)\,ds(\zeta) + o(\epsilon).
  \end{equation*}
  
  \item In Laplacian growth dynamics, a domain containing the origin grows with outward velocity given by
  $V(\zeta) = \partial g/\partial n(\zeta, 0)$. In these circumstances, an infinitesimal time step $dt$ causes a domain variation with
  $\epsilon p = V\,dt$. As such, the variation formula
  yields
  \begin{equation*}
    \frac{dg}{dt}(z,0) = -\int_{\partial D} \frac{\partial g}{\partial n}(z,\zeta)\left[\frac{\partial g}{\partial n}(0,\zeta)\right]^2  \,ds(\zeta).
  \end{equation*}
  A generalization of Laplacian growth dynamics is elliptic growth dynamics, wherein we use the Green function for $L = \nabla \lambda \nabla$ and take
  $p(\zeta) = \lambda(\zeta)\partial g/\partial n(\zeta, 0)$. In this case the previous formula is modified into
  \begin{equation*}
    \frac{dg}{dt}(z,0) = -\int_{\partial D} \frac{\partial g}{\partial n}(z,\zeta)\left[\lambda(\zeta)\frac{\partial g}{\partial n}(0,\zeta)\right]^2  \,ds(\zeta).
  \end{equation*}
  
  \item There is an alternative way of defining the Green function, wherein $g(z,w) = -\ln|z-w| + O(1)$ as $z\to w$. Taking this as a definition leads to a few minor changes. Firstly,
  $\Delta g_w = -2\pi \delta_w$ and the Poisson-Jensen formula reads
  \begin{equation*}
    u(z) = \int_{\partial D} u(\zeta)\, d\omega(z,\zeta) -  \frac{1}{2\pi}\int_D g(z,\xi) \Delta u(\xi)\, dA(\xi).
  \end{equation*}
  This leads to the variation formula
  \begin{equation*}
    g^*(z,w) - g(z,w) = \frac{\epsilon}{2\pi}\int_{\partial D} \frac{\partial g}{\partial n}(z,\zeta)\frac{\partial g}{\partial n}(w,\zeta) p(\zeta)\,ds(\zeta) + o(\epsilon).
  \end{equation*}
\end{enumerate}

\section{Operator Variation of the Green Function}

Now we turn to situations wherein the domain is fixed but rather the underlying operator is somehow close to the Laplacian. We will make repeated use
of the following integral operator.

\begin{definition}
  Given a domain $D$ define the operator $T$ to be integration against the Green function of $D$:
  \begin{equation*}
    T\phi(z) = \int_D \phi(\xi)g(z,\xi)\, dA(\xi).
  \end{equation*}
\end{definition}

\begin{lemma}
  The operator $T$ is a bounded linear map from $L^2(D)$ into $C(D)$, the space of continuous functions on $D$.
\end{lemma}
\begin{proof}
  Notice that $T$ is a left inverse to the Laplacian. By elliptic regularity, $T$ maps $L^2$ into $H^2$.
  But in the plane, Sobolev imbedding implies that $H^2(D) \subset C(D)$. 
\end{proof}

\subsection{Helmholtz}

Suppose we'd like to approximate the Green function $g^*$ of the Helmholtz operator $H = \Delta - a$, where $a$ is a small constant. In terms of the original
Green function $g$ we can derive the following.
\begin{theorem}\label{Helmholtz}
  Let $D$ be a bounded domain in $\C$ and fix $w\in D$.
  The Green function $g^*$ of the Helmholtz operator $\Delta-a$ satisfies
  \begin{equation*}
     g^*_w - g_w = aTg_w + o(a)
  \end{equation*}
  as $a\to 0$, where the convergence of $o(a)$ is uniform in $z$ for each fixed $w$. Furthermore, a full series expansion is given by
  \begin{equation*}
    g^*_w = \sum_{n=0}^\infty a^n T^n g_w.
  \end{equation*}
\end{theorem}
\begin{proof}
  The function $g^*_w$ solves the boundary value problem
  \begin{equation*}
    \begin{cases}
    \Delta g^*_w &= \delta_w + ag^*_w \quad \textrm{ in } D\\
    g^*_w &= 0 \quad \textrm{ on }\partial D
    \end{cases}
  \end{equation*}
  By definition of the classical Green function, this implies that
  \begin{equation*}
    g^*_w(z) = T(\delta_w(\xi) + ag^*_w(\xi)) = g_w(z) + aTg^*_w(\xi).
  \end{equation*}
  If we define the error function $f = g^*_w-g_w$ the previous equation becomes
  \begin{equation}
    f = aTg^*_w = aTg_w + aTf. \label{helmEqn}
  \end{equation}
  It remains to show that $Tf \to 0$ uniformly as $a\to 0$.
  For all $a$ smaller than $\|T\|^{-1}$,
  \begin{equation}
    f = a(I-aT)^{-1}Tg_w. \label{HelmholtzExact}
  \end{equation}
  Since $g_w\in L^2$ so is $(I-aT)^{-1}Tg_w$ and we can take $a\to 0$ to find
  \begin{equation*}
    \|f\|_2 \leq |a| \|(I-aT)^{-1}Tg_w\|_2 \to 0.
  \end{equation*}
  That is, $f\to 0$ in $L^2$ as $a\to 0$. The lemma implies that $Tf \to 0$ in $C(D)$, hence uniformly.
  In addition we obtain the full series expansion
  \begin{equation*}
    g^*_w = \sum_{n=0}^\infty a^n T^n g_w
  \end{equation*}
  if we use the Neumann series for $(I-aT)^{-1}$ in \eqref{HelmholtzExact}.
\end{proof}

\subsection{Schr\"odinger}
In the previous section we deduced a perturbation formula for a Schr\"odinger operator with constant potential function.
With slight adjustments to the proof we can handle a general potential function; the operator $\Delta - u$ also possesses a perturbation formula
when $u$ is small.
\begin{theorem}\label{Schrodinger}
  Let $D$ be a bounded domain in $\C$ and fix $w\in D$. Suppose that $p$ is a smooth scalar function defined in a neighborhood of $D$
  with corresponding multiplication operator $P$.
  The Green function $g^*$ of the Schr\"odinger operator $\Delta-\epsilon p$ satisfies
  \begin{equation*}
     g^*_w - g_w = \epsilon TPg_w + o(\epsilon)
  \end{equation*}
  as $\epsilon \to 0$, where the convergence of $o(\epsilon)$ is uniform in $z$ for each fixed $w$. Furthermore, a full series expansion is given by
  \begin{equation*}
    g^*_w = \sum_{n=0}^\infty \epsilon^n (TP)^n g_w.
  \end{equation*}
\end{theorem}
\begin{proof}
  The first steps of the proof are similar to the Helmholtz case.
  The function $g^*_w$ solves the boundary value problem
  \begin{equation*}
    \begin{cases}
    \Delta g^*_w &= \delta_w + \epsilon pg^*_w \quad \textrm{ in } D\\
    g^*_w &= 0 \quad \textrm{ on }\partial D
    \end{cases}
  \end{equation*}
  By definition of the classical Green function, this implies that
  \begin{equation}
    g^*_w = T(\delta_w + \epsilon p g^*_w) = g_w + \epsilon TPg^*_w = g_w + \epsilon TPg_w + \epsilon TP(g^*_w-g_w). \label{SchrodEqn}
  \end{equation}
  Define the error function $f = g^*_w-g_w$, so
  \begin{equation*}
    f = \epsilon TPg_w + \epsilon TPf.
  \end{equation*}
  It remains to show that $TPf \to 0$ uniformly as $\epsilon\to 0$.
  
  As in our treatment of the Helmholtz operator, $T$ is a bounded linear map $L^2 \to H^2$. Since $D$ is bounded, $p\in L^\infty(D)$; thus
  $P$ maps $L^2\to L^2$. Together this means that $TP$ continuously maps $L^2 \to H^2$ and the rest of proof will proceed analogously as before. For all
  $\epsilon$ sufficiently small we find that
  \begin{equation}
    f = \epsilon (I-\epsilon TP)^{-1}TPg_w. \label{SchrodExact}
  \end{equation}
  Since $g_w\in L^2$ we can take $\epsilon\to 0$:
  \begin{equation*}
    \|f\|_2 \leq |\epsilon| \|(I-\epsilon TP)^{-1} TP g_w\|_2 \to 0.
  \end{equation*}
  Thus $f\to 0$ in $L^2$ and $TPf \to 0$ in $H^2 \subset C(D)$. We conclude that $Tf \to 0$ uniformly.
  The full series expansion again follows from expanding $(1-\epsilon TP)^{-1}$ in a Neumann series.
\end{proof}

\subsection{Laplace--Beltrami}
Next we consider the problem of finding the Green function on a domain for the Laplace--Beltrami operator $L = \nabla \lambda \nabla$, where $\lambda$ is a smooth positive function.
If $\lambda$ is close to unity that the Green function $g^*$ can be approximated by the Green function $g$ of the Laplacian.
The following result makes this notion precise.
\begin{theorem}\label{Beltrami}
  Fix $w\in D$ and suppose that $p$ is a smooth scalar function in a neighborhood of $D$.
  Given $\epsilon > 0$ we can define $\lambda(z) = 1+\epsilon p(z)$. Then as $\epsilon \to 0$
  the Green function $g^*$ for $L=\nabla \lambda \nabla$ satisfies
  \begin{equation}
    g^*(z,w) - g(z,w) = \epsilon \int_D p(\xi) \nabla g(z,\xi) \cdot \nabla g(\xi,w)\, dA(\xi) + o(\epsilon), \label{Beltrami1}
  \end{equation}
  where all derivatives are with respect to $\xi$. Furthermore, the error term converges uniformly for each fixed $w$.
  An alternate formula is also true:
  \begin{equation}
    g^*(z,w) - g(z,w) = -\epsilon g(z,w)\left(\frac{p(z) + p(w)}{2}\right) + \frac{\epsilon}{2}\int_D g_z g_w \Delta p \, dA + o(\epsilon). \label{Beltrami2}
  \end{equation}
\end{theorem}
\begin{proof}
  We first give a heuristic argument. Once again, if we begin with the correct integration by parts formula we can formally manipulate our way to the
  perturbation formula above.
  From the Green identity
  \begin{equation*}
    \int_{\partial D} \lambda\left(u\frac{\partial v}{\partial n} - v\frac{\partial u}{\partial n}\right)\, ds = \int_{D} \left(uL v - vL u\right)\, dA
  \end{equation*}
  we let $u=g_z$ and $v=g^*_w$ to get
  \begin{equation*}
    \int_D g_z Lg^*_w\, dA = \int_D g^*_w Lg_z\, dA.
  \end{equation*}
  Note that $Lg^*_w = \delta_w$ and $Lg_z = \lambda\delta_z + \nabla\lambda\cdot\nabla g_z$ to find
  \begin{equation*}
    g(z,w) = \lambda(z)g^*(z,w) + \int_D g^*_w \nabla\lambda \cdot\nabla g_z\, dA,
  \end{equation*}
  whence
  \begin{align}
    g^*(z,w) - g(z,w) &= -\epsilon p(z)g^*(z,w) - \epsilon\int_D g^*_w \nabla p\cdot\nabla g_z \, dA \nonumber \\ 
    &= -\epsilon p(z)g(z,w) - \epsilon\int_D g_w \nabla p \cdot\nabla g_z \, dA + \epsilon R_\epsilon, \label{lbpf}
  \end{align}
  where we have defined
  \begin{equation*}
    R_\epsilon = -p(z)(g^*(z,w)-g(z,w)) - \int_D (g^*_w - g_w) \nabla p \cdot\nabla g_z \, dA.
  \end{equation*}
  Here we end the rigor and claim that $R_\epsilon=o(1)$. Notice that this is the same gap that appeared in the derivation
  of Hadamard's formula above. Using the claim \eqref{lbpf} becomes
  \begin{equation}
    g^*(z,w) - g(z,w) = -\epsilon p(z)g(z,w) - \epsilon\int_D g_w \nabla p \cdot\nabla g_z \, dA + o(\epsilon). \label{lbpf2}
  \end{equation}
  We can simplify this via integration by parts. Note that $g(\xi,w) = 0$ for $\xi\in\partial D$; therefore
  \begin{align*}
    \int_D \nabla p \cdot(g_w \nabla g_z )\, dA
    &=
     - \int_D p \nabla\cdot (g_w\nabla g_z)\, dA \\
    &= - \int_D p g_w\Delta g_z \, dA -  \int_D p \nabla g_z \cdot \nabla g_w \, dA\\
    &= - p(z) g(z,w)  -  \int_D p \nabla g_z\cdot \nabla g_w\, dA.
  \end{align*}
  Inserting this into \eqref{lbpf2} gives our first perturbation formula \eqref{Beltrami1}.
  To deduce \eqref{Beltrami2} we take \eqref{lbpf2} and reverse the roles of $z$ and $w$:
  \begin{equation*}
    g^*(z,w) - g(z,w) = -\epsilon p(w)g(z,w) - \epsilon\int_D g_z \nabla p \cdot\nabla g_w \, dA + o(\epsilon).
  \end{equation*}
  Adding this half of this equation to half of \eqref{lbpf2}, we recognize a product rule.
  \begin{equation*}
    g^*(z,w) - g(z,w) = -\epsilon g(z,w)\left(\frac{p(z)+p(w)}{2}\right) - \frac{\epsilon}{2}\int_D \nabla p \cdot\nabla (g_z g_w)\, dA + o(\epsilon).
  \end{equation*}
  Since $g(z,\xi)g(\xi,w) = 0$ on $\partial D$, integrating by parts one more time gives the result.
  
  We now proceed to derive \eqref{Beltrami2} in another fashion, sidestepping the need to consider the error term $R_\epsilon$. Essentially
  we will use a change of variables to relate the problem to that of the Sch\"odinger operator.
  
  Define the functions $G^*_w = g^*_w\sqrt{\lambda}$ and $u = \lambda^{-1/2}\Delta(\lambda^{1/2})$. Notice that $G^*_w = 0$ on $\partial D$ and
  \begin{equation}
    \nabla \lambda = \nabla(\sqrt{\lambda}\sqrt{\lambda}) = 2\sqrt{\lambda}\nabla\sqrt{\lambda}. \label{RootEqn}
  \end{equation}
  With this in mind we compute:
  \begin{align*}
    (\Delta - u)G^*_w &= \Delta(g^*_w\sqrt{\lambda}) - ug^*_w\sqrt{\lambda} \\
    &= g^*_w\Delta\sqrt{\lambda} + \sqrt{\lambda}\Delta g^*_w + 2\nabla g^*_w\cdot\nabla \sqrt{\lambda} - g^*_w\Delta\sqrt{\lambda} \\
    &= \frac{\lambda\Delta g^*_w + \nabla\lambda\cdot\nabla g^*_w}{\sqrt{\lambda}} \\
    &= \delta_w /\sqrt{\lambda} \\
    &= \delta_w /\sqrt{\lambda(w)}.
  \end{align*}
  From this we conclude that $\sqrt{\lambda(w)}G^*_w$ is the Green function of $D$ for the operator $\Delta - u$. We can use our perturbation formula
  for this operator if we can determine how $u$ depends upon $\epsilon$. We begin with a computation; from equation \eqref{RootEqn} we find that
  \begin{equation*}
    |\nabla \sqrt{\lambda}|^2 = \left|\frac{\nabla\lambda}{2\sqrt{\lambda}}\right|^2 = \frac{|\nabla\lambda|^2}{4\lambda}.
  \end{equation*}
  Now we have what we need to relate $u$ and $\lambda$ without roots.
  \begin{equation*}
    \Delta\lambda = \Delta \sqrt{\lambda}\sqrt{\lambda} = 2\sqrt{\lambda}\Delta\sqrt{\lambda} + 2|\nabla\sqrt{\lambda}|^2 = 2\lambda u + \frac{|\nabla\lambda|^2}{2\lambda},
  \end{equation*}
  from which it follows that
  \begin{equation*}
    u = \frac{\Delta \lambda}{2\lambda} - \frac{|\nabla\lambda|^2}{4\lambda^2} = \frac{\epsilon\Delta p/2}{1+\epsilon p} - \frac{\epsilon^2|\nabla p|^2/4}{(1+\epsilon p)^2}.
  \end{equation*}
  From here we expand in series to find
  \begin{align*}
    u &= \frac{\epsilon\Delta p}{2}\sum_{n=0}^\infty (-1)^n (\epsilon p)^n - \frac{\epsilon^2|\nabla p|^2}{4}\sum_{n=0}^\infty n(-1)^n(\epsilon p)^{n-1} \\
    &= \sum_{n=0}^\infty (-1)^n \left(\frac{\epsilon \Delta p}{2} + \frac{n\epsilon |\nabla p|^2}{4p}\right) (\epsilon p)^n \\
    &= \frac{\epsilon \Delta p}{2} + o(\epsilon),
  \end{align*}
  where the error term is uniformly convergent for all small $\epsilon$. Therefore we have
  \begin{equation*}
    \sqrt{\lambda(w)\lambda(z)}g^*(z,w) = g(z,w) + \int_D \left(\frac{\epsilon \Delta p}{2} + o(\epsilon)\right)g_z g_w\, dA + o(\epsilon).
  \end{equation*}
  As a map $L^\infty \to L^\infty$ the integral operator
  \begin{equation*}
    \phi \mapsto \int_D \phi g_z g_w\, dA
  \end{equation*}
  is continuous, so we deduce that
  \begin{equation*}
    \int_D o(\epsilon) g_z g_w\, dA = o(\epsilon)
  \end{equation*}
  uniformly. We conclude that
  \begin{equation*}
    \sqrt{\lambda(w)\lambda(z)}g^*(z,w) = g(z,w) + \frac{\epsilon}{2}\int_D g_z g_w \Delta p \, dA + o(\epsilon).
  \end{equation*}
  The last step requires a first order approximation of $(\lambda(w)\lambda(z))^{-1/2}$ in $\epsilon$. We find that
  \begin{align*}
    \frac{1}{\sqrt{\lambda(w)\lambda(z)}} &= (1+\epsilon p(w))^{-1/2} (1+ \epsilon p(z))^{-1/2} \\
    &= \left(1 - \frac{\epsilon p(z)}{2} + o(\epsilon)\right)\left(1 - \frac{\epsilon p(w)}{2} + o(\epsilon)\right) \\
    &= 1-\epsilon\left(\frac{p(z)+p(w)}{2}\right) + o(\epsilon).
  \end{align*}
  Therefore
  \begin{align*}
    g^*(z,w) &= \left[1-\epsilon\left(\frac{p(z)+p(w)}{2}\right) + o(\epsilon)\right]\left[g(z,w) + \frac{\epsilon}{2}\int_D g_z g_w \Delta p\, dA + o(\epsilon)\right] \\
    &= g(z,w) - \epsilon g(z,w)\left(\frac{p(z)+p(w)}{2}\right) + \frac{\epsilon}{2}\int_D g_z g_w \Delta p\, dA + o(\epsilon),
  \end{align*}
  as desired.
\end{proof}

\begin{example}
Again consider $\D$, the unit disk in $\C$. If we take $\lambda(\xi) = 1+\epsilon |\xi|^2$ then
\begin{equation*}
  \delta g(z,w) = \Re\frac{1}{4\pi^2}
    \int_\D |\xi|^2\left(\frac{1}{\xi-w} + \frac{\bar{w}}{1-\bar{w}\xi}\right)\cdot \overline{\left(\frac{1}{\xi-z} + \frac{\bar{z}}{1-\bar{z}\xi}\right)}\, dA(\xi).
\end{equation*}
We can evaluate the integral as follows.
First we assume that $z,w\neq 0$ and notice that
\begin{equation*}
  \xi\left(\frac{1}{\xi-w} + \frac{\bar{w}}{1-\bar{w}\xi}\right) = \frac{w}{\xi-w} + \frac{1}{1-\bar{w}\xi},
\end{equation*}
whence the integrand can be rewritten as
\begin{equation}
    \frac{w\bar{z}}{(\xi-w)(\bar{\xi}-\bar{z})}  + \frac{w}{(1-z\bar{\xi})(\xi-w)} + \frac{\bar{z}}{(1-\xi\bar{w})(\bar{\xi}-\bar{z})} + \frac{1}{(1-\xi\bar{w})(1-\bar{\xi}z)}. \label{integrand}
\end{equation}
We will integrate each of these four terms in turn. Following \cite{ExpTrans}, we recognize the exponential transform of $\D$:
\begin{align*}
  \int_{\D} \frac{w\bar{z}}{(\xi-w)(\bar{\xi}-\bar{z})}\, dA(\xi) &= -\pi w\bar{z} \log E_\D(w,z) \\
  &= \pi w\bar{z}\log(1-w\bar{z}) - 2\pi w\bar{z} \log|z-w|.
\end{align*}
Here we have used the fact that $\log(ab) = \log a + \log b \pmod {i\R}$; ultimately we only desire the real part of each integral.
This handles the first term in \eqref{integrand}. Next consider
\begin{equation*}
   \int_{R} \frac{w}{(1-z\bar{\xi})(\xi-w)} \, dA(\xi),
\end{equation*}
where $R = \D\setminus D(w,\delta)$ and $\delta$ is a small positive number. Since $(1-z\bar{\xi})^{-1}$ is antianayltic throughout $\D$ we can write
\begin{align*}
   \int_{R} \frac{w}{(1-z\bar{\xi})(\xi-w)} \, dA(\xi) &= -\frac{w}{z}\int_{R} \dbar\left( \frac{\log(1-z\bar{\xi})}{\xi-w}\right) \, dA(\xi) \\
   &= -\frac{w}{z}\int_{\partial R} \frac{\log(1-z\bar{\xi})}{\xi-w} \, \frac{d\xi}{2i}.
\end{align*}
Considering orientation of the inner circle, we have
\begin{equation*}
   -\frac{w}{z}\int_{\partial R} \frac{\log(1-z\bar{\xi})}{\xi-w} \, \frac{d\xi}{2i}
   = -\frac{w}{z}\int_{\partial \D} \frac{\log(1-z\bar{\xi})}{\xi-w} \, \frac{d\xi}{2i} +\frac{w}{z}\int_{\partial D(w,\delta)} \frac{\log(1-z\bar{\xi})}{\xi-w} \, \frac{d\xi}{2i}
\end{equation*}
To evaluate the first of these integrals we remark that $z/\xi\in \D$ for $\xi\in \partial \D$ so a series expansion gives
\begin{align*}
  -\frac{w}{z}\int_{\partial \D} \frac{\log(1-z/\xi)}{\xi-w} \, \frac{d\xi}{2i} &= \frac{\pi w}{z}\sum_{k=1}^\infty \int_{\partial \D} \frac{z^k}{k \xi^k(\xi - w)}\, \frac{d\xi}{2\pi i} \\
  &= \frac{\pi w}{z}\sum_{j,k=1}^\infty \int_{\partial \D} \frac{z^k w^{j-1}}{k \xi^{k+j}}\, \frac{d\xi}{2\pi i} \\
  &= 0.
\end{align*}
The dominated convergence theorem implies that
\begin{align*}
  \lim_{\delta\to 0} \frac{w}{z}\int_{\partial D(w,\delta)} \frac{\log(1-z\bar{\xi})}{\xi-w} \, \frac{d\xi}{2i} &=
  \lim_{\delta \to 0}\frac{w}{z}\int_0^{2\pi} \log(1-z\bar{w} - \delta z e^{-i\theta}) \, \frac{d\theta}{2} \\
  &= \frac{\pi w}{z}\log(1-z\bar{w}).
\end{align*}
Thus
\begin{equation*}
   \int_{\D} \frac{w}{(1-z\bar{\xi})(\xi-w)} \, dA(\xi) = \frac{\pi w}{z}\log(1-z\bar{w}).
\end{equation*}
Similarly,
\begin{equation*}
   \int_{\D} \frac{\bar{z}}{(1-\xi\bar{w})(\bar{\xi}-\bar{z})} \, dA(\xi) = \frac{\pi \bar{z}}{\bar{w}}\log(1-z\bar{w}).
\end{equation*}
The last term in \eqref{integrand} can be evaluated as such:
\begin{align*}
   \int_{\D} \frac{1}{(1-\xi\bar{w})(1-\bar{\xi}z)} \, dA(\xi) &= -\frac{1}{z} \int_{\D} \dbar\left(\frac{\log(1-z\bar{\xi})}{1-\xi \bar{w}}\right)\, dA(\xi) \\
   &= -\frac{1}{z} \int_{\partial\D} \frac{\log(1-z/\xi)}{1-\xi \bar{w}}\frac{d\xi}{2i} \\
   &= \sum_{j,k=1}^\infty \int_{\partial \D} \frac{z^{k-1} \bar{w}^{j-1}}{k\xi^{k-j+1}} \, \frac{d\xi}{2i} \\
   &= \pi\sum_{k=1}^\infty \frac{(z\bar{w})^{k-1}}{k} \\
   &= -\frac{\pi}{z\bar{w}}\log(1-z\bar{w}).
\end{align*}
Altogether we have
\begin{equation}
  \delta g(z,w) = \frac{1}{4\pi}\Re\left[\left(z\bar{w}+\frac{\bar{z}}{\bar{w}} + \frac{w}{z} - \frac{1}{z\bar{w}}\right)\log(1-z\bar{w}) - 2z\bar{w}\log|z-w|\right],
  \label{BeltramiAnswer}
\end{equation}
where we have taken complex conjugates of various terms to combine them better. Note that this formula only depends upon rotationally-invariant quantities such as $z\bar{w}$, as can be expected.

Next we assume that $w=0$ and compute
\begin{align*}
  \delta g(z,0) &= \Re\frac{1}{4\pi^2}
    \int_\D |\xi|^2 \bar{\xi}^{-1}\cdot \left(\frac{1}{\xi-z} + \frac{\bar{z}}{1-\bar{z}\xi}\right)\, dA(\xi) \\
    &= \Re \frac{1}{4\pi^2}\int_\D \left(\frac{\xi}{\xi-z} + \frac{\bar{z}\xi}{1-\bar{z}\xi}\right)\, dA(\xi).
\end{align*}
The second term is a harmonic function in $\xi$; therefore we have
\begin{equation*}
    \int_\D \frac{\bar{z}\xi}{1-\bar{z}\xi} \, dA(\xi) = \pi \frac{\bar{z}\xi}{1-\bar{z}\xi}\bigg|_{\xi=0} = 0.
\end{equation*}
Given a small $\delta >0$ define $D_\delta = \D\setminus D(z,\delta)$ and use Stokes' theorem to write
\begin{align*}
  \int_{D_\delta} \frac{\xi}{\xi - z}\, dA(\xi) &= \int_{D_\delta} \dbar\left(\frac{|\xi|^2}{\xi-z}\right)\, dA(\xi) \\
  &= \int_{\partial \D} \frac{|\xi|^2}{\xi - z}\frac{d\xi}{2i} - \int_{\partial D(z,\delta)} \frac{|\xi|^2}{\xi - z}\frac{d\xi}{2i} \\
  &= \pi - \frac{1}{2}\int_0^{2\pi} |z+\delta e^{i\theta}|^2\, d\theta \\
  &\to \pi(1-|z|^2)
\end{align*}
as $\delta\to 0$. From this we conclude that
\begin{equation*}
  \delta g(z,0) = \frac{1-|z|^2}{4\pi}.
\end{equation*}
We remark the same result is obtained by taking $w\to 0$ in equation \eqref{BeltramiAnswer}.
\end{example}

\bibliographystyle{amsplain}

\begin{thebibliography}{99}



\bibitem{Garabedian} P. R. Garabedian, \textsl{Partial Differential Equations}, Chelsea, New York 1986.

\bibitem{ExpTrans} B. Gustafsson, M. Putinar, \textsl{An exponential transform and regularity of free boundaries in two dimensions}, Ann. 
Scuola Norm. Sup. Pisa Cl. Sci. (4), \textbf{26} (1998), 507--543.

\bibitem{Hadamard} J. Hadamard, \textsl{M\'emoire sur le probl\`eme d'analyse relatif \`a l'\'equilibre des plaques \'elastiques encastr\'ees},
M\'emoires present\'es par divers savants \`a l'Acad\'emie des Sciences \textbf{33} (1908), 1--128.

\bibitem{SS} E. Schippers, W. Staubach, \textsl{Variation of Neumann and Green functions under homotopies of the boundary}, Israel J. Math. \textbf{173} (2009), 279--303.

\end{thebibliography}

		 \textsc{Department of Mathematics,
		 University of California, Santa Barbara,
		 CA 93106} \\
\textit{E-mail:} \texttt{cmart07@math.ucsb.edu}

\end{document}